\documentclass[a4paper,12pt,reqno]{amsart}
\usepackage{amsmath,amssymb}
\usepackage{url}

\usepackage{setspace}

\newtheorem{Def}{Definition}
\newtheorem{Rmk}{Remark}

\newtheorem{Thm}{Theorem}

\newtheorem{Lem}{Lemma}

\newcommand{\Q}[1]{\mathbb{Q}(\sqrt{#1})}
\newcommand{\Z}{\mathbb{Z}}

\newcommand{\frakO}{\mathfrak{O}}

\newcommand{\im}{\operatorname{Im}}

\newcommand{\qf}[1]{\left\langle #1 \right\rangle}

\newcommand{\rank}{\operatorname{rank}}
\newcommand{\gen}{\operatorname{gen}}

\newcommand{\s}[1]{\mathfrak{s}#1}
\newcommand{\n}[1]{\mathfrak{n}#1}
\newcommand{\vol}[1]{\mathfrak{v}#1}
\newcommand{\binlattice}[4]{\begin{pmatrix}
  #1 & #2 \\
  #3 & #4
\end{pmatrix}}
\newcommand{\terlattice}[9]{\begin{pmatrix}
  #1 & #2 & #3 \\
  #4 & #5 & #6 \\
  #7 & #8 & #9
\end{pmatrix}}
\newcommand{\terdiagonal}[3]{\begin{pmatrix}
  #1 & 0 & 0 \\
  0 & #2 & 0 \\
  0 & 0 & #3
\end{pmatrix}}

\newcommand{\conj}[1]{\overline{#1}}
\newcommand{\comega}{\conj\omega}
\newcommand{\nto}{\nrightarrow}
\newcommand{\nleftto}{\hspace{.5mm}\not{\hspace{-1.5mm}\leftarrow}}
\newcommand{\nequiv}{\not\equiv}

\title[2-Universal Hermitian Lattices]{2-universal Hermitian lattices \\ over imaginary quadratic fields}

\author[Myung-Hwan Kim]{Myung-Hwan Kim}
\address{Myung-Hwan Kim\\
Department of Mathematical Sciences, Seoul National University, San 56-1 Shinrim-dong Gwanak-gu, Seoul,
151-747, Korea}%
\email{mhkim@math.snu.ac.kr}

\author[Poo-Sung Park]{Poo-Sung Park}
\address{Poo-Sung Park\\
School of Computational Sciences,
Korea Institute institute for Advanced Study, Hoegiro 87, Dongdaemun-gu, Seoul, 130-722, Korea}%
\email{sung@kias.re.kr}


\thanks{The authors were partially supported by KRF(2005-070-c00004)}

\begin{document}

\begin{abstract}
A positive definite Hermitian lattice is said to be $2$-universal if
it represents all positive definite binary Hermitian lattices. We
find all $2$-universal ternary and quaternary Hermitian lattices
over imaginary quadratic number fields.
\end{abstract}

\maketitle

\section{Introduction}
We call a positive definite integral quadratic form \emph{universal} if it represents all positive integers. Then Lagrange's Four Square Theorem means that the sum of four squares is universal. In 1930, Mordell \cite{Mordell} generalized this notion to a \emph{$2$-universal} quadratic form: a positive definite integral quadratic form that represents all binary positive definite integral quadratic forms, and showed that the sum of five squares is $2$-universal. In this direction, we refer the readers to \cite{K} and \cite{KKO1,KKO2}.

As another generalization of universal quadratic forms, universal Hermitian forms have been studied. This was initiated by Earnest and Khosravani. They defined a \emph{universal Hermitian} form as the one representing all positive integers, and found 13 universal binary Hermitian forms over imaginary quadratic fields of class number 1 \cite{E-K}. The list of binary universal Hermitian forms has been completed by Iwabuchi \cite{Iwabuchi}, Jae-Heon Kim and the second author \cite{K-P}. The simple and unified proofs was recently obtained by the second author \cite{Park}. In this paper, we study $2$-universal Hermitian forms. We prove that there are finitely many $2$-universal ternary and quaternary Hermitian forms over imaginary quadratic fields, and find them all (sections \ref{sec:2UniversalTernary} and \ref{sec:2UniversalQuaternary}).

A notable recent progress in the representation theory of quadratic forms is the so called \emph{Fifteen Theorem} of Conway-Schneeberger \cite{Conway}, which states: \emph{a positive definite quadratic form is universal if it represents positive integers up to 15}. This fascinating result was improved by Bhargava \cite{Bhargava}, who proved analogies for other infinite subsets of positive integers like the set of all primes, the set of all positive odd integers and so on. Kim et al. \cite{KKO1,KKO2} recently proved the finiteness theorem for representability and provided a $2$-universal analogy of the Fifteen Theorem. Recently Kim, Kim and the second author proved the Fifteen Theorem for universal Hermitian lattices. In section \ref{sec:FinitenessTheorem}, we obtain a criterion for $2$-universality of Hermitian forms over several imaginary quadratic fields.


\section{Notations and Symbols}

Let $E=\Q{-m}$ for a square-free positive integer $m$ and $\frakO = \frakO_E = \mathbb Z[\omega]$ be its ring of integers, where $\omega=\omega_m=\sqrt{-m}$ \,or $\frac{1+\sqrt{-m}}2$ if $m \equiv 1, 2$ or $3$ $\pmod{4}$, respectively.

For a prime $p$ we define $E_p := E \otimes_\mathbb{Q} \mathbb{Q}_p$. Then the ring $\frakO_p$ of integers of $E_p$ is defined $\frakO \otimes_\Z \Z_p$. If $p$ is inert or ramifies in $E$, then $E_p = \mathbb{Q}_p(\sqrt{-m})$ and $\alpha\otimes\beta = \alpha\beta$ with $\alpha \in E$ and $\beta \in \mathbb{Q}_p$. If $p$ splits in $E$, then $E_p = \mathbb{Q}_p \times \mathbb{Q}_p$ and $\alpha\otimes\beta = (\alpha\beta,\conj{\alpha}\beta)$ where $\conj{\cdot}$ is the canonical involution. Thus $E_p$ allows the unique involution $\conj{\alpha\otimes\beta} = \conj{\alpha}\otimes\beta$ \cite{Gerstein}.

\begin{Def}
Let $F=E$ or $E_p$. A \emph{Hermitian space} is a finite-dimensional vector space $V$ over $F$ equipped with a sesqui-linear map $H: V\times V \to F$ satisfying the following conditions:
\begin{enumerate}
    \item $H(\mathbf{x},\mathbf{y}) = \conj{H(\mathbf{y},\mathbf{x})}$,
    \item $H(a\mathbf{x},\mathbf{y}) = aH(\mathbf{x},\mathbf{y})$,
    \item $H(\mathbf{x}_1+\mathbf{x}_2,\mathbf{y}) = H(\mathbf{x}_1,\mathbf{y}) + H(\mathbf{x}_2,\mathbf{y})$,
\end{enumerate}
We simply denote $H(\mathbf{v},\mathbf{v})$ by $H(\mathbf{v})$ and call it the \emph{(Hermitian) norm} of $\mathbf{v}$.
\end{Def}

From (1)$\sim$(3) follow %
\[
(2') \ H(\mathbf{x},b\mathbf{y})=\conj{b}H(\mathbf{x},\mathbf{y}) \mbox{ and } (3') \  H(\mathbf{x},\mathbf{y}_1+\mathbf{y}_2) = H(\mathbf{x},\mathbf{y}_1)+H(\mathbf{x},\mathbf{y}_2).
\]

\begin{Def}
Let $R=\frakO$ or $\frakO_p$. A \emph{Hermitian $R$-lattice} is an $R$-module $L$ equipped with a sesqui-linear map $H$ such that $H(L,L)\subseteq R$.
\end{Def}

If $L$ is free with a basis $\{\mathbf{v}_1, \dotsc, \mathbf{v}_n\}$, then we define
\[
	M_L:=(H(\mathbf{v}_i,\mathbf{v}_j))_{n\times n}\,,
\] 
and call it the \emph{Gram matrix} of $L$. We often identify $M_L$ with the lattice $L$. If $M_L$ is diagonal, we simply write $L = \qf{a_1, \dotsc, a_n}$, where $a_i = H(\mathbf{v}_i)$ for $i=1,2,\dotsc,n$. The determinant of $M_L$ is called the \emph{discriminant} of $L$, denoted by $dL$. The {\it discriminant} $dV$ of a Hermitian space $V$ over $F$ is defined analogously as the determinant of $M_V$ and is well-defined up to $N(\dot{F})$, where $N$ is the norm map on $F$ defined by $N(a)=a\conj{a}$ for $a\in F$. If $dL$ is a unit, we call $L$ \emph{unimodular}. By $FL$, we mean the Hermitian space $V=F\otimes_R L$ where $L$ is nested. We define the \emph{rank} of $L$ by $\rank L:=\dim_F FL$. By the \emph{scale} $\s{L}$ and the \emph{norm} $\n{L}$ of $L$, we mean the $R$-modules generated by the subsets $H(L,L)$ and $H(L)$, respectively.

It is well known that an $R$-lattice $L$ can be written as
\begin{equation}
    L = \mathfrak{a}_1\mathbf{v}_1 + \dotsb +     \mathfrak{a}_n\mathbf{v}_n
\label{ideal_expression1}
\end{equation}%
for vectors $\mathbf{v}_1, \dotsc, \mathbf{v}_n\in L$ and ideals $\mathfrak{a}_1,\dotsc,\mathfrak{a}_n \subseteq R$. We define the \emph{volume} of $L$ by
\[
    \vol{L}: = \mathfrak{a}_1\conj{\mathfrak{a}}_1 \cdots \mathfrak{a}_n\conj{\mathfrak{a}}_n \det(H(\mathbf{v}_i,\mathbf{v}_j)).
\]
We call $L$ \emph{modular} if $\vol{L} = (\s{L})^n$. The expression (\ref{ideal_expression1}) of an $\frakO$-lattice $L$, which is not necessarily free, can be transformed into the form
\begin{equation}
    L = \frakO\mathbf{w}_1 +  \dotsb + \frakO\mathbf{w}_{n-1} + \mathfrak{a}\mathbf{w}_n
\label{ideal_expression2}
\end{equation}%
for some vectors $\mathbf{w}_1, \dotsc, \mathbf{w}_n\in L$ and an ideal $\mathfrak{a}\subseteq \frakO$ \cite[81:5]{OM}. If $L$ is not free, or equivalently, if $\mathfrak{a}$ is not principal, then $\mathfrak{a}$ is generated by two elements, say $\alpha, \beta \in \frakO$. Therefore (\ref{ideal_expression2}) may be rewritten as
\[
    L = \frakO\mathbf{w}_1 + \frakO\mathbf{w}_2 + \dotsb + \frakO\mathbf{w}_{n-1} + \frakO \alpha
   \mathbf{w}_n + \frakO \beta\mathbf{w}_n.
\]
We may treat $L$ as if it were a free lattice with basis $\{\mathbf{w}_1, \dots, \mathbf{w}_{n-1}, \alpha\mathbf{w}_n, \beta\mathbf{w}_n\}$. The rank of $L$, however, is still $n$ not $n+1$. In this case the formal Gram matrix of $L$ is defined as
\[
    M_L =
    \begin{pmatrix}
    H(\mathbf{v}_1,\mathbf{v}_1)        & \dotsc & H(\mathbf{v}_1,\alpha\mathbf{v}_n)        & H(\mathbf{v}_1,\beta\mathbf{v}_n) \\
    \vdots            & \ddots & \vdots                   & \vdots \\
    H(\alpha\mathbf{v}_n,\mathbf{v}_1) & \dotsc & H(\alpha\mathbf{v}_n,\alpha\mathbf{v}_n) & H(\alpha\mathbf{v}_n,\beta\mathbf{v}_n) \\
    H(\beta\mathbf{v}_n,\mathbf{v}_1)  & \dotsc & H(\beta\mathbf{v}_n,\alpha\mathbf{v}_n)  & H(\beta\mathbf{v}_n,\beta\mathbf{v}_n)
    \end{pmatrix}.
\]

Let $\ell$ and $L$ be two (free or nonfree) Hermitian $R$-lattices whose (formal) Gram matrices are $M_\ell \in M_{m \times m}(R)$ and $M_L \in M_{n \times n}(R)$ respectively. We say that $L$ \emph{represents} $\ell$, denoted by $\ell \to L$, if there exists a suitable $X \in M_{m\times n}(R)$ such that $M_\ell = XM_LX^\ast$, where $X^\ast$ is the conjugate transpose of $X$. The two lattices are said to be \emph{isometric}, denoted by $\ell \cong L$, if they represent each other. For Hermitian spaces $v$ and $V$, $v \to V$ and $v \cong V$ are defined analogously.


\section{Escalation Method}

A Hermitian $\frakO$-lattice $L$ is called \emph{positive $($definite$)$} if $H(v)>0$ for all $v\in L, \,v\ne 0$. From here on, we assume that every Hermitian $\frakO$-lattice is positive unless stated otherwise.

\begin{Def}
A positive Hermitian $\frakO$-lattice $L$ is called \emph{$2$-universal} if it represents all binary positive Hermitian $\frakO$-lattices.
\end{Def}

Let $\ell$ be a binary Hermitian $\frakO$-lattice and let $G$ be a minimal generating set of vectors of $\ell$. We let $S(G)$ be the largest norm of vectors in $G$ and define
\[
	S(\ell):=\min_G \{\,S(G)\,:\,G \mbox{ is a minimal generating set of } \ell \,\}.
\]

\begin{Def}
Let $L$ be a Hermitian $\frakO$-lattice that is not $2$-universal. We define the \emph{truant} of $L$ by
\[
	T(L):=\min_{\ell} \{\,S(\ell)\,:\,\ell \mbox{ is binary and } \ell\nrightarrow L\,\}.
\]
\end{Def}

For example, let $E=\Q{-1}$ and $L = \qf{1,1}$ be a Hermitian $\frakO$-lattice. Then $\ell_1=\qf{1,3}\nrightarrow L$ and $S(\ell_1) = 3$. But 3 is not the truant of $L$ because $S(\ell_2)=2$ for $\ell_2 = \binlattice2112\nrightarrow L$. It is clear that 2 is the truant of $L$.

\begin{Def}
An \emph{escalation} of a non-$2$-universal lattice $L$ is defined to be any lattice which is generated by $L$ and a vector whose norm is equal to $T(L)$. An \emph{escalation lattice} is a lattice which can be obtained as the result of a sequence of successive escalations from the zero-dimensional lattice.
\end{Def}

We now construct $2$-universal lattices using the escalation method starting from the zero-dimensional lattice.

Let $L_0$ denote the zero dimensional lattice. Since $T(L_0)=1$, the first escalation lattice is $L_1=\qf{1}$. Since $\qf{1,1} \nrightarrow L_1$, we again have $T(L_1)=1$ and hence the second escalation lattice is $L_2=\qf{1,1}$.  

In order to find the third escalation lattices, we consider $\ell=\binlattice2112$. Suppose $\ell\to L_2$, that is,
\begin{align*}
    \binlattice2112 &=
    \binlattice{a}{b}{c}{d}
    \binlattice1001
    \binlattice%
        {\conj{a}}{\conj{c}}
        {\conj{b}}{\conj{d}}
    =
    \binlattice%
        {a\conj{a}+b\conj{b}}
        {a\conj{c}+b\conj{d}}
        {\conj{a}c+\conj{b}d}
        {c\conj{c}+d\conj{d}}.
\end{align*}
for some $a,b,c,d\in \frakO$. From this we obtain
\[
    a\conj{a}+b\conj{b} = 2, \ \ c\conj{c}+d\conj{d} = 2     \mbox{ and } a\conj{c}+b\conj{d} = 1,
\]
which imply that $a\conj{a} = b\conj{b} =c\conj{c} = d\conj{d} = 1$. Elements of norm 1 in $\frakO$ are
\[
\left\{\begin{array}{ll} \pm1&\quad\mbox{if}\quad m\ne 1,3\\
\pm1, \pm \sqrt{-1}&\quad\mbox{if}\quad m= 1\\
\pm1, \pm\frac{1\pm\sqrt{-3}}2&\quad\mbox{if}\quad m=3.
\end{array}\right.
\]
So, if $m \ne 3$, then $a\conj{c}+b\conj{d}=1$ cannot be satisfied and thus $\binlattice2112 \nrightarrow L_2$. If $m=3$, then $\qf{1,2}\nrightarrow L_2$ because $2$ is not a norm of any element in $\Z[\frac{1+\sqrt{-3}}{2}]$. Therefore, we may conclude that $T(L_2)=2$ for any $m$ and hence that the third escalation lattices are\,:
\[
    L_{3;1}=\qf{1,1,1}, \ L_{3;2}=\qf{1,1,2}
\]
if free, and \ $L_{3;3}=\qf{1,1}\perp\binlattice{2}{b}{\conj{b}}{c} \mbox{ \ with \ } 2c-b\conj{b}=0$ \, if not free.


\section{2-Universal Ternary Hermitian Lattices}
\label{sec:2UniversalTernary}

Let's assume that $m\ne 1,2,3,7,11$ for the time being. Observe that there is no element of norm 2 and 3 in $\frakO$ because $\omega\comega \ge 4$ under this assumption. Using this observation, it is easy to show that
\[
	\qf{1,3}\nrightarrow L_{3;1}=\qf{1,1,1} \ \ 
	\mbox{ and } \ \
	\qf{3,3}\nrightarrow L_{3;2}=\qf{1,1,2}.
\]

Consider the nonfree case\,: $L_{3;3}=\qf{1,1}\perp\binlattice{2}{b}{\conj{b}}{c} \mbox{ \ with \ } 2c-b\conj{b}=0$. This case may occur only when $h_E\ne 1$, where $h_E$ is the ideal class number of $E$. So, $m\ne 1,2,3,7,11,19,43,67,163$, necessarily. By a simple reduction, we may assume that $b=\omega$ or $-1+\omega$. Suppose now that $L_{3;3}$ is $2$-universal.\\

\noindent Subcase 1) Let $b=\omega$. Then the corresponding Hermitian form is
\begin{align*}
    H_{\omega}(x,y,z,u) &= x\conj{x} + y\conj{y} + 2z\conj{z} + \omega z\conj{u} +     \comega\conj{z}u + cu\conj{u} \\
    &= x\conj{x} + y\conj{y} + \frac12(2z+\comega u)(2\conj{z}+\omega     \conj{u}).
\end{align*}
Since $\qf{1,5} \to L_{3;3}$, the equation $y\conj{y} + \frac12(2z+\comega u)(2\conj{z}+\omega\conj{u})=5$ should be solvable over $\frakO$. The following table lists those $m$'s that satisfy the equation and $2c-\omega\comega=0$ solvable over $\frakO$.\\

\begin{center}
\begin{tabular}{ccc}
    $y\conj{y}$ & $(2z+\comega u)(2\conj{z}+\omega\conj{u})$ & $m$ \\
    \hline
    0 & 10 & 6,10,15,31,39 \\
    1 & 8 & 23,31 \\
    \hline
\end{tabular}
\end{center}\vskip 0.2cm%
For each of $m$ in the table, it is easy to verify that
$\binlattice2112 \nrightarrow L_{3;3}$.

\noindent Subcase 2) Let $b=-1+\omega$. Then the corresponding
Hermitian form is
\begin{align*}
    H_{-1+\omega}(x,y,z,u) &= x\conj{x} + y\conj{y} + 2z\conj{z}     + (-1+\omega) z\conj{u} +     (-1+\comega)\conj{z}u + cu\conj{u} \\
    &= x\conj{x} + y\conj{y} + \frac12[2z+(-1+\comega)     u][2\conj{z}+(-1+\omega)\conj{u}]
\end{align*}
Since $\qf{1,5} \to L_{3;3}$, the equation $y\conj{y} + \frac12[2z+(-1+\comega)u][2\conj{z}+(-1+\omega)\conj{u}]=5$ should be solvable over $\frakO$. The following table lists those $m$'s that satisfy the equation and $2c-(-1+\omega)(-1+\comega)=0$ solvable over
$\frakO$.\\

\begin{center}
\begin{tabular}{ccc}
    $y\conj{y}$ & $[2z+(-1+\comega)u][2\conj{z}+(-1+\omega)\conj{u}]$ & $m$ \\
    \hline
    0 & 10 & 15,31,39 \\
    1 & 8 & 23,31 \\
    5 & 0 & 5 \\
    \hline
\end{tabular}
\end{center}\vskip 0.2cm%
For each of $m$ in the table, it is easy to verify that $\binlattice2112 \nrightarrow L_{3;3}$.\\

In summary, we proved\,:
\begin{Thm}
There is no $2$-universal ternary Hermitian $\frakO$-lattice over the imaginary quadratic fields $\Q{-m}$ if $m \ne 1,2,3,7,11$. Moreover, every $2$-universal ternary Hermitian $\frakO$-lattice over $\Q{-m}$ is free.
\end{Thm}

We now assume that $m=1,2,3,7,11$. We'll give new names for convenience to the remaining candidates for $2$-universal ternary Hermitian $\frakO$-lattices as follows\,: $$I:=L_{3;1}=\qf{1,1,1} \quad \mbox{ and }\quad J:=L_{3;2}=\qf{1,1,2}.$$ We eliminate more candidates by finding binary lattices that cannot be represented as follows.\\

\begin{center}
\begin{tabular}{ll}
$m=1$\,: &
    $\binlattice2112 \nrightarrow J$\,;\\
$m=2$\,: &
    $\binlattice2{-1+\omega}{-1+\comega}2\nrightarrow I, J$\,; \\
$m=7$\,: &
    $\binlattice2112\nrightarrow J$\,; \\
$m=11$\,: &
    $\binlattice2{\omega}{\comega}2, \binlattice2{-1+\omega}{-1+\comega}2\nrightarrow I, J$\,.\\
\end{tabular}
\end{center}
\vskip 0.2cm

The remaining candidates for $2$-universal ternary Hermitian
$\frakO$-lattices are
$$\begin{array}{ll}
    I=\qf{1,1,1} & \text{ over } \Q{-1}\,; \\
    I=\qf{1,1,1}, \ J=\qf{1,1,2} & \text{ over } \Q{-3}\,; \\
    I=\qf{1,1,1} & \text{ over } \Q{-7}
\end{array}$$
and they are indeed $2$-universal.

\begin{Thm}\label{Thm:111_over_Q_sqrt_negative_1}
The ternary Hermitian $\frakO$-lattice $I=\qf{1,1,1}$ over $\Q{-1}$ is $2$-universal.
\end{Thm}

\begin{proof}
Let $\ell$ be a binary Hermitian $\frakO$-lattice. Since the class number of $I$ is $1$ \cite{Otremba}, it is enough to show that $\ell_p \to L_p$ for every prime $p$.

If $p$ is split, then $\ell_p \to I_p$ by \cite[1.8]{Gerstein}. If $E_p/\mathbb{Q}_p$ is an unramified quadratic extension, $\ell_p \to L_p$ by \cite[Theorem 4.4]{Johnson}. If $E_p/\mathbb{Q}_p$ is a non-dyadic ramified quadratic extension, $\ell_p \to L_p$ by \cite[Theorem 5.5]{Johnson}. If $p=2$, $\ell_p \to L_p$ by \cite[Theorem 9.4]{Johnson}. The representation at the archimedean prime spot is clear. Therefore, $I=\qf{1,1,1}$ is $2$-universal.
\end{proof}

\begin{Thm}\label{Thm:111,112_over_Q_sqrt_negative_3} 
    The ternary Hermitian $\frakO$-lattices $I=\qf{1,1,1}$ and     $J=\qf{1,1,2}$ over $\Q{-3}$ are $2$-universal.
\end{Thm}

\begin{proof}
Let $E=\Q{-3}$. Since $h(I)=h(J)=1$ \cite{Otremba}, it is enough to prove the local universality for every prime $p$. The local universality of $I$ and the local universality of $J$ at each prime spot are checked by \cite[1.8]{Gerstein} if $p$ is split, \cite[Theorem 4.4]{Johnson} if $E_p/\mathbb{Q}_p$ is an unramified quadratic extension, or \cite[Theorem 5.5]{Johnson} if $E_p/\mathbb{Q}_p$ is a ramified quadratic extension.
\end{proof}

\begin{Thm}\label{Thm:111_over_Q_sqrt_negative_7}
    The ternary Hermitian $\frakO$-lattice $I = \qf{1,1,1}$ over $\Q{-7}$     is $2$-universal.
\end{Thm}

\begin{proof}
Let $E=\Q{-7}$ and let $\ell$ be a binary Hermitian $\frakO$-lattice. We obtain the local representation $\ell_p \rightarrow I_p$ at every $p$ by \cite[1.8]{Gerstein}, \cite[Theorem 4.4]{Johnson}, or \cite[Theorem 5.5]{Johnson}. But, we cannot say that $\ell \to I$ because $h(I) = 2$ \cite{Schiemann}.

We have the other class in $\gen{I}$, say the class of $T$, where
\[
    T := \terlattice%
        {2}{1}{\omega}
        {1}{2}{1}
        {\comega}{1}{2}.
\]
Thus, in order to prove the $2$-universality of $I$, we should prove that $\ell\rightarrow T$ implies $\ell\rightarrow I$. To this end, let's assume $\ell \to T$, that is,
\[
    \ell =
    \begin{pmatrix}
        a_1 & a_2 & a_3 \\
        b_1 & b_2 & b_3
    \end{pmatrix}
    T
    \begin{pmatrix}
        \conj{a}_1 & \conj{b}_1 \\
        \conj{a}_2 & \conj{b}_2 \\
        \conj{a}_3 & \conj{b}_3
    \end{pmatrix}.
\]
for some $a_1,a_2,a_3,b_1,b_2,b_3\in \frakO$. Let $X:=\begin{pmatrix} a_1 & a_2 & a_3 \\ b_1 & b_2 & b_3 \end{pmatrix}$.

We call an algebraic integer in $\Q{-7}$ \emph{even-type} if it is divisible by $\omega$, and \emph{odd-type} otherwise. It is easy to show that two algebraic integers $a$ and $b$ are of the same type if and only if the sum $a+b$ is of even-type. For a column $\begin{pmatrix} a\\b \end{pmatrix}$ of $X$, call it \emph{EE-type} if both $a$ and $b$ are of even-type, \emph{EO-type} if $a$ is of even-type and $b$ is of odd-type, \emph{OE-type} if $a$ is of odd-type and $b$ is of even-type, and \emph{OO-type} if both $a$ and $b$ are of odd-type.

Firstly, assume that the first column of $X$ is of EE-type. By Letting $a_1=a_1'\omega$ and $b_1=b_1'\omega$, we get
\[
    \ell = X \terlattice
        {2}{1}{\omega} {1}{2}{1} {\comega}{1}{2} X^* =
    \begin{pmatrix}
        a_1' & a_2 & a_3 \\
        b_1' & b_2 & b_3
    \end{pmatrix}
    \terlattice
        {4}{\omega}{-2+\omega}
        {\comega}{2}{1}
        {-2+\comega}{1}{2}
   \begin{pmatrix}
        \conj{a}_1' & \conj{b}_1' \\
        \conj{a}_2 & \conj{b}_2 \\
        \conj{a}_3 & \conj{b}_3
    \end{pmatrix}
\]
and
\[
    T' :=
    \terlattice%
        {1 }{1}{\comega} {0 }{1}{-1     } {-1}{0}{-1     }
    I
    \terlattice%
        {1     }{0 }{-1} {1     }{1 }{0 } {\omega}{-1}{-1}
    =
    \terlattice
        {4}{\omega}{-2+\omega}
        {\comega}{2}{1}
        {-2+\comega}{1}{2},
\]
which implies that $\ell\rightarrow I$. Similarly, if the second or
third column is of EE-type, $T$ can be replaced by
\[
    \terlattice
        {2       }{\comega}{\omega}
        {\omega}{4      }{\omega}
        {\comega }{\comega}{2     }
    \text{ or }
    \terlattice
        {2}{1     }{2      }
        {1}{2     }{\comega}
        {2}{\omega}{4      }
\]
as $T'$, respectively, and both are represented by $I$.

Secondly, assume that none of the columns of $X$ is of EE-type. If
the first two columns of $X$ are of the same type, then we take
\[
    T' =
    \terlattice%
        {6        }{3\comega}{1+\omega}
        {3\omega  }{4       }{\omega  }
        {1+\comega}{\comega }{2       }
\]
using $a_2 = (a_2-a_1)+a_1 = a_2'\omega+a_1$ and $b_2 = (b_2-b_1)+b_1 = b_2'\omega+b_1$. Similarly, if the last two columns of $X$ are of the same type, then we take
\[
    T' =
    \terlattice%
        {2}{1+\omega}{2}
        {1+\comega}{6}{3\comega}
        {2}{3\omega}{4},
\]
and if the first and the third columns of $X$ are of the same type, we take
\[
    T' =
    \terlattice%
        {5        }{2     }{2+2\comega}
        {2        }{2     }{\comega   }
        {2+2\omega}{\omega}{4         }.
\]
All these $T'$s are represented by $I$.

Finally, assume that no two columns of $X$ are of the same type. Then all three types other than the EE-type should occur in clumns of $X$. In this case, we take
\[
    T' =
    \terlattice%
        {5        }{4+\omega}{2+2\comega}
        {4+\comega}{6       }{3\comega  }
        {2+2\omega}{3\omega }{4         }
\]
using $a_3 = (a_3-a_1-a_2)+a_1+a_2 = a_3'\omega+a_1+a_2$ and $b_3 = (b_3-b_1-b_2)+b_1+b_2 = b_3'\omega+b_1+b_2$. This is also can be represented by $I$.

Therefore, we may conclude that if $\ell\rightarrow T$, then $\ell\rightarrow I$, as desired.
\end{proof}


\section{$2$-universal quaternary Hermitian lattices}
\label{sec:2UniversalQuaternary}

Recall that the escalated ternary lattices are
\[
    L_{3;1} = \qf{1,1,1} \quad\text{ and }\quad L_{3;2} = \qf{1,1,2}.
\]
Now we escalate these ternary lattices to construct $2$-universal lattices. We already found $2$-universal lattices $\qf{1,1,1}$ over $\Q{-1}$, $\Q{-3}$, and $\Q{-7}$. Thus we may assume that $m \ne 1,3,7$ in escalating $\qf{1,1,1}$ and $m \ne 3$ in escalating
$\qf{1,1,2}$.

Since $\omega\comega \ge 3$, $L=\qf{1,1,1,a}$ does not represent $\binlattice{2}\omega\comega{c}$ with $\frac{\omega\comega}2 < c < \omega\comega$. All $\Q{-m}$ except for $m=2$ are excluded, since $\omega\comega \ge 3$ for $m \ne 1,2,3,7$. When $m=2$, we verify that
\[
    \binlattice{2}{-1+\omega}{-1+\comega}2 \nto L
\]
for any $a$. Thus there are no universal lattices escalated from $\qf{1,1,1}$.

Consider $\qf{1,1,2}$. We have that $\qf{1,1,2}$ is universal over $\Q{-3}$ \cite{K-P}. So assume that $m \ne 3$. Then since $\binlattice{2}{1}{1}{2} \nto \qf{1,1,2}$, the truant is $2$. Thus the next escalation lattice is
\[
    L_4 = \qf{1,1} \perp \binlattice{2}b{\conj{b}}2
\]
up to isometry. If $b=0$, then $\binlattice{2}{1}{1}{2} \nto \qf{1,1,2,2}$ for $m \ne 3$. Hence we may assume that $b \ne 0$. We obtain an inequality $1 \le b\conj{b} < 4$ for $L$ to be positive definite. Then the feasible lattices are as follows:
\begin{center}
\begin{tabular}{ll}
    $m$ & $b$ \\ \hline
    1,2,7,11 & 1, $\omega_{m}$, $-1+\omega_{m}$ \\
    otherwise & 1 \\ \hline
\end{tabular}
\end{center}

When $m=5,6,10$ or $m \ge 13$, $\omega\comega \ge 4$ and thus
\[
    \binlattice{2}{\omega}{\comega}c \nto
    \qf{1,1}\perp\binlattice{2}{1}{1}{2}
\]
if $\frac{\omega\comega}2 < c < \omega\comega$.

So (new) $2$-universal quaternary Hermitian lattices are of the form $\qf{1,1}\perp\binlattice{2}{\ast}{\ast}{2}$ and $m=1,2,7,11$. We find a binary lattice which cannot be represented by each quaternary lattice.

\begin{center}
\begin{tabular}{ll}
$\Q{-1}$: &
    $\qf{1,1}\perp\binlattice{2}{-1+i}{-1-i}2 \nleftto
        \binlattice{2}{1}{1}{2}, \binlattice{2}\omega{\comega}2$;\\
$\Q{-2}$: &
    $\qf{1,1}\perp\binlattice{2}{1}{1}{2} \nleftto
        \binlattice{2}{-1+\omega_2}{-1+\comega_2}2$, \\
\null &
    $\qf{1,1}\perp\binlattice{2}{\omega_2}{\comega_2}{2} \nleftto
        \binlattice{2}{1}{1}{2}, \binlattice{2}{-1+\omega_2}{-1+\comega_2}2$; \\
$\Q{-7}$: &
    $\qf{1,1}\perp\binlattice{2}{1}{1}{2} \nleftto
        \binlattice3{1+\omega_7}{1+\comega_7}3$, \\
\null &
    $\qf{1,1}\perp\binlattice{2}{\omega_7}{\comega_7}2 \nleftto
        \binlattice{2}{1}{1}{2}$; \\
$\Q{-11}$: &
    $\qf{1,1}\perp\binlattice{2}{1}{1}{2} \nleftto
        \binlattice{2}{\omega_{11}}{\comega_{11}}{2}$
\end{tabular}
\end{center}

Now let us investigate nonfree quaternary lattices. Note that nonfree lattices are possible when the ideal class number of imaginary quadratic field $\Q{-m}$ is bigger than $1$ and thus $m=5,6,10$ or $m \ge 13$. Then $\omega_{m}\comega_{m} \ge 4$.

Assume a quaternary lattice $L$ is nonfree and $2$-universal. We have that
\[
    L = \qf{1,1,1} \perp \binlattice{\ast}{\ast}{\ast} {\ast} \text{
    or } L = \qf{1,1} \perp \terlattice{2}{\ast}{\ast}
    {\ast}{\ast}{\ast} {\ast}{\ast}{\ast}
\]
and one of principal minors vanishes for the formal Gram matrix $M_L$ to be positive semi-definite.

Consider the first case. Let $L = \qf{1,1,1}\perp L_0$. If $\min L_0 \ge 3$, then $L$ cannot represent $\binlattice{2}113$. Thus we can write
\[
    L = \qf{1,1,1}\perp\binlattice{2}\omega\comega{c} \text{ or } L =
    \qf{1,1,1}\perp\binlattice{2}{-1+\omega}{-1+\comega}c
\]
with $\det M_L=0$ up to isometry. We can write these lattices as
\[
    L = \qf{1,1,1}\perp(2,b)\qf{\frac12}
    = \frakO \mathbf{v}_1 + \frakO \mathbf{v}_2
        + \frakO \mathbf{v}_3 + (2,b)\frakO \mathbf{v}_4
\]
with $b = \comega$ or $-1+\comega$. $L$ must represent the binary lattice $\binlattice3003$. Then $3 \to \qf{1,1,1}$ and $3 \to (2,b)\qf{\frac12}$ since $\alpha\conj\alpha \ge 4$ unless $\alpha = \pm1$ . Note that the vector $\mathbf{v}_1 + 2\mathbf{v}_4$ also make norm $3$ but it is not orthogonal to any other vector of norm 3. Since $6 \in N((2,b)\frakO)$ is possible only when $m = 6,15,23$, the candidates are as following. But none of them are $2$-universal:

\begin{center}
\begin{tabular}{ll}
    $\Q{-6}$: &
    $\qf{1,1,1}\perp\binlattice{2}{\omega_{6}}{\comega_{6}}3
    \nleftto \binlattice{2}{-1+\omega_{6}}{-1+\comega_{6}}3$; \\
    $\Q{-15}$: &
    $\qf{1,1,1}\perp\binlattice{2}{\omega_{15}}{\comega_{15}}2
    \nleftto \binlattice4{-1+2\omega_{15}}{-1+2\comega_{15}}4$,
    \\
    \null &
    $\qf{1,1,1}\perp\binlattice{2}{-1+\omega_{15}}{-1+\comega_{15}}2
    \nleftto \binlattice4{-1+2\omega_{15}}{-1+2\comega_{15}}4$;
    \\
    $\Q{-23}$: &
    $\qf{1,1,1}\perp\binlattice{2}{\omega_{23}}{\comega_{23}}3
    \nleftto \binlattice3{1+\omega_{23}}{1+\comega_{23}}3$,
    \\
    \null &
    $\qf{1,1,1}\perp\binlattice{2}{-1+\omega_{23}}{-1+\comega_{23}}3
    \nleftto \binlattice3{1+\omega_{23}}{1+\comega_{23}}3$
\end{tabular}
\end{center}

Now consider the second case: $L =
\qf{1,1}\perp\terlattice{2}{\ast}{\ast}{\ast}{\ast}{\ast}
{\ast}{\ast}{\ast}$. Let us consider $\ell =
\binlattice{2}{1}{1}{2}$. Since $\ell \nto \qf{1,1,2}$, its truant
is $2$. Thus, the escalation lattice is
\[
    \qf{1,1} \perp \binlattice{2}{1}{1}{2}.
\]

Next, the lattice does not represent a binary lattice $\ell = \binlattice{2}\omega\comega{c}$, where $c$ is the smallest integer satisfying $2c-\omega\comega > 0$. That is, $c$ is the truant and thus
\[
    L = \qf{1,1} \perp \terlattice21\omega 12b \comega{\conj{b}}{c}.
\]

One of leading principal minors should be 0 for $L$ to be positive semi-definite. Thus
\begin{align*}
    \det M_L &= 3c-2b\conj{b}+\omega\conj{b}+\comega b -
    2\omega\comega\\
    &= c + (2c-\omega\comega) - b\conj{b} -
    (b-\omega)(\conj{b}-\comega) \\
    &= 0
\end{align*}
and the fact that $2c-\omega\comega = 1$ or $2$ yields
\begin{equation}\label{Eqn:inequality}
    b\conj{b} + (b-\omega)(\conj{b}-\comega)
    =
    \begin{cases}
    \displaystyle c + 1 = \frac{\omega\comega}2+\frac32 & \text{ if $\omega\comega$ is
    odd,} \\[2ex]
    \displaystyle c + 2 = \frac{\omega\comega}2+3 & \text{ if $\omega\comega$ is
    even.}
    \end{cases}
\end{equation}
Note that $b$ or $b-\omega$ has nonzero $\omega$-part. Choose that number and write it as $s+t\omega$ with $t\ne0$. Then
\begin{align*}
    b\conj{b} + (b-\omega)(\conj{b}-\comega)
    &\ge (s+t\omega)(s+t\comega) \\
    &= s^2 + st(\omega+\comega) + t^2\omega\comega \\
    &\ge \omega\comega.
\end{align*}
Thus the equality (\ref{Eqn:inequality}) can be satisfied when $\omega\comega = 4,6$. The following lattices are all candidates up to isometry, but none of them is $2$-universal:

\begin{center}
\begin{tabular}{ll}
    $\Q{-6}$: &
    $\qf{1,1}\perp
    \terlattice{2}{2}{4} {1}{\omega_6}{0} {1}{\comega_6}{0}
    \nleftto \binlattice{2}{0}{0}{3}$; \\
    $\Q{-15}$: &
    $\qf{1,1}\perp
    \terlattice{2}{2}{3} {1}{\omega_{15}}{1} {1}{\comega_{15}}{1}
    \nleftto \binlattice{4}{-1+2\omega_{15}}{-1+2\comega_{15}}{4}$; \\
    $\Q{-23}$: &
    $\qf{1,1}\perp
    \terlattice{2}{2}{4} {1}{\omega_{23}}{0} {1}{\comega_{23}}{0}
    \nleftto \binlattice{2}{0}{0}{3}$
\end{tabular}
\end{center}

In summary candidates for $2$-universal quaternary Hermitian $\frakO$-lattices are
\[
\begin{array}{ll}
    \qf{1,1}\perp\binlattice{2}{1}{1}{2} & \text{ over } \Q{-1}\,; \\
    \qf{1,1}\perp\binlattice{2}{-1+\omega_2}{-1+\comega_2}{2} & \text{ over } \Q{-2}\,; \\
    \qf{1,1}\perp\binlattice{2}{\omega_{11}}{\comega_{11}}{2} & \text{ over } \Q{-11}
\end{array}
\]
and they are indeed $2$-universal.


\begin{Thm}\label{quaternary_lattice_for_m_equal_to_1}
    The Hermitian lattice $\qf{1,1}\perp\binlattice{2}{1}{1}{2}$ over     $\Q{-1}$ is $2$-universal.
\end{Thm}
To prove this we need some setups. Let $E=\Q{-1}$. If an algebraic integer $\alpha\in\Q{-1}$ is divisible by $1+i$, we     call $\alpha$ \emph{even type}. If not, we call it \emph{odd type}. If we gather all vectors of even norm in a Hermitian lattice $L$ over $E$, they form a lattice, say $L_{e}$. We construct a basis of $L_e$ using this facts.

\begin{Lem}\label{Lem:basis_of_Le}
    Let $\{\mathbf{v}_1, \mathbf{v}_2, \dotsb, \mathbf{v}_n \}$ be a basis of $L$. Assume that     $H_L(\mathbf{v}_j)$ is odd for each $j$. Then the lattice $L_e$ of even     vectors has a basis as follows:
    \begin{align*}
        \mathbf{v}_1' &= \mathbf{v}_1 - \mathbf{v}_2, \\
        \mathbf{v}_2' &= \mathbf{v}_2 - \mathbf{v}_3, \\
         & \dotsb \\
        \mathbf{v}_{n-1}' &= \mathbf{v}_{n-1} - \mathbf{v}_n \\
        \mathbf{v}_n' &= \mathbf{v}_{n-1} + \sqrt{-1} \mathbf{v}_n
    \end{align*}
\end{Lem}
\begin{proof}
Let $\mathbf{v} = a_1 \mathbf{v}_1 + \dotsb + a_n \mathbf{v}_n$ be a vector in $L_e$. Then the sum of all $a_j$'s is even type since $H_L(\mathbf{v}_j)$ are all odd. Thus
\begin{align*}
    \mathbf{v} &= a_1 \mathbf{v}_1 + \dotsb + a_n \mathbf{v}_n \\
    &= a_1 (\mathbf{v}_1 - \mathbf{v}_2) + (a_1 + a_2) (\mathbf{v}_2 - \mathbf{v}_3) \\
    &+ \dotsb \\
    &+ (a_1 + a_2 + \dotsb + a_{n-1}) (\mathbf{v}_{n-2} - \mathbf{v}_{n-1}) \\
    &+ \left(\frac{i}{1+i}(a_1 + a_2 +
    \cdots + a_{n-1} + a_n)-a_n\right) (\mathbf{v}_{n-1} - \mathbf{v}_n) \\
    &+ \frac{1}{1+i}(a_1 + a_2 +
    \cdots + a_{n-1} + a_n) (\mathbf{v}_{n-1} + i \mathbf{v}_n).
\end{align*}
\end{proof}

Now we prove the $2$-universality.

\begin{proof}[Proof of Theorem \ref{quaternary_lattice_for_m_equal_to_1}]
Let $L = \qf{1,1,2}$. Its class number is $1$ by \cite{Otremba}. Let $\ell$ be a binary Hermitian $\mathfrak{O}$-lattice. If $p \ne 2$, $\ell_p \to L_p$ by \cite[1.8]{Gerstein} or \cite[Theorem4.4]{Johnson}. 

Let $p=2$. Then $\ell_p \to L_p$ unless $\ell_p$ is unimodular and $\n{\ell_p} = \frakO_p$ by \cite[Theorem 9.4]{Johnson}. If $\ell_p$ is unimodular and $\n{\ell_p} = \frakO_p$, then $\ell \cong \binlattice{2a}{b}{\conj{b}}{2c}$ with $2 \nmid b\conj{b}$. To see such $\ell$ is representable, we consider the sublattice $\qf{1}\perp\binlattice{2}{1}{1}{2}$. Ternary lattices of discriminant $3$ are all decomposable \cite{Zhu}. Thus other lattice is the only $\qf{1,1,3}$ and the two lattices compose a genus.

It is easy to check that $\ell_p \to \qf{1,1,3}_p$ for every $p$. Thus $\ell$ is represented by $\qf{1}\perp\binlattice{2}{1}{1}{2}$ or $\qf{1,1,3}$. Suppose $\ell \to L:=\qf{1,1,3}$. Note that all vectors contained in $\ell$ have even norms. Thus we may assume that $\ell \to L_e$ instead of $\ell \to L$.

Let $\{\mathbf{v}_1, \mathbf{v}_2, \mathbf{v}_3\}$ be a basis of $L$ such that
\[
    H_L(\mathbf{v}_1) = 1, H_L(\mathbf{v}_2) = 1, H_L(\mathbf{v}_3) = 3, H_L(\mathbf{v}_j, \mathbf{v}_k) = 0     \text{ for } j \ne k.
\]
Then using Lemma \ref{Lem:basis_of_Le}, we obtain that
\begin{align*}
    &H_L(\mathbf{v}_1') = 2, H_L(\mathbf{v}_2') = 4, H_L(\mathbf{v}_3') = 4, \\
    &H_L(\mathbf{v}_1', \mathbf{v}_2') = -1, H_L(\mathbf{v}_1', \mathbf{v}_3') = -1, H_L(\mathbf{v}_2', \mathbf{v}_3') = 1+i,
\end{align*}
and thus we can write
\[
    L_e \cong \terlattice%
        {2 }{-1 }{-1 }
        {-1}{4  }{1+i}
        {-1}{1-i}{4  }.
\]

Hence we conclude that $\ell \to \qf{1,1}\perp\binlattice2112$ since $L_e \to \qf{1,1}\perp\binlattice2112$ from
\[
    \terlattice%
        {2 }{-1 }{-1 }
        {-1}{4  }{1+i}
        {-1}{1-i}{4  }
    =
    \begin{pmatrix}
        0 & 0 & 0 & 1 \\
        0 & 1+i & -1 & 0 \\
        1 & -i & -1 & 0
    \end{pmatrix}
    \begin{pmatrix}
        1 & 0 & 0 & 0 \\
        0 & 1 & 0 & 0 \\
        0 & 0 & 2 & 1 \\
        0 & 0 & 1 & 2
    \end{pmatrix}
    \begin{pmatrix}
        0 & 0 & 1 \\
        0 & 1-i & i \\
        0 & -1 & -1 \\
        1 & 0 & 0
    \end{pmatrix}.
\]
\end{proof}

\begin{Thm}
    The Hermitian lattice
    $\qf{1,1}\perp\binlattice{2}{-1+\sqrt{-2}}{-1-\sqrt{-2}}{2}$
    over $\Q{-2}$ is $2$-universal.
\end{Thm}
\begin{proof}
Consider a lattice $\qf{1,1,2}$. Then any positive binary Hermitian lattice can be represented by one of $\gen(\qf{1,1,2})$. Since genus of the integral quadratic form $x_1^2+x_2^2+y_1^2+y_2^2+2z_1^2+2z_2^2$ induced by the Hermitian lattice has a minimum $\le 2$ from the Hermite constant $\gamma_6 = \frac{2}{\sqrt[6]{3}}$. So genus of Hermitian lattice $\qf{1,1,2}$ has a minimum $1$ or $2$. By using this fact we can find the genus which consists of four classes:
\begin{align*}
    &\terlattice{1}{0}{0} {0}{1}{0} {0}{0}{2},
    \terlattice{2}{0}{\omega} {0}{2}{-1+\omega}
    {\comega}{-1+\comega}{3},\\
    &\terlattice{2}{\omega}{-1+\omega} {\comega}{2}{-1+\omega} {-1+\comega}{-1+\comega}{7},
    \terlattice{2}{-1}{-1+\omega} {-1}{5}{-1+2\omega} {-1+\comega}{-1+2\comega}{5}.
\end{align*}
These lattices are all represented by the quaternary lattice. So it is $2$-universal.
\end{proof}

\begin{Thm}
    The Hermitian lattice     $\qf{1,1}\perp\binlattice{2}{\omega}{\comega}{2}$ over $\Q{-11}$     is $2$-universal, where $\omega = \frac{1+\sqrt{-11}}2$.
\end{Thm}
\begin{proof}
It is clear that $\qf{1,1,1}_p$ represents any binary lattice locally for every $p$. Since the genus of $\qf{1,1,1}$ consists of
\[
    \terdiagonal111 \text{ and }
    \terlattice10002{\omega}0{\comega}2.
\]
by \cite{Schiemann}, it is enough to show that a binary lattice represented by $\qf{1,1,1}$ can also be represented by $\qf{1,1}\perp\binlattice{2}{\omega}{\comega}{2}$.

We call a number \emph{$0$-type} if it is divisible by $\omega$. If $n \in \pm1+\Z\omega$, we call it \emph{$1$-type} or \emph{$-1$-type}, respectively. Then the difference of same type numbers is 0-type, that is, the difference is divisible by $\omega$. If a number is $1$-type and the other is $-1$-type, we say that they have \emph{opposite types}.

Suppose that a binary lattice $\ell$ can be represented as
\[
    \ell =
    \begin{pmatrix}
        a_1 & a_2 & a_3 \\
        b_1 & b_2 & b_3
    \end{pmatrix}
    \terdiagonal111
    \begin{pmatrix}
        \conj{a}_1 & \conj{b}_1 \\
        \conj{a}_2 & \conj{b}_2 \\
        \conj{a}_3 & \conj{b}_3
    \end{pmatrix}.
\]
If $a_1$ and $b_1$ are both $0$-type, replace them by $\omega a_1'$ and $\omega b_1'$, respectively. The matrix calculation shows us that
\[
    \ell =
    \begin{pmatrix}
        a_1' & a_2 & a_3 \\
        b_1' & b_2 & b_3
    \end{pmatrix}
    \terlattice
        300
        010
        001
   \begin{pmatrix}
        \conj{a}_1' & \conj{b}_1' \\
        \conj{a}_2 & \conj{b}_2 \\
        \conj{a}_3 & \conj{b}_3
    \end{pmatrix}
\]
and
\[
    \terdiagonal311 =
    \begin{pmatrix}
        0 & 0 & 1 & -1 \\
        1 & 0 & 0 & 0 \\
        0 & 1 & 0 & 0
    \end{pmatrix}
    \begin{pmatrix}
        1 & 0 & 0 & 0 \\
        0 & 1 & 0 & 0 \\
        0 & 0 & 2 & \omega \\
        0 & 0 & \comega & 2
    \end{pmatrix}
    \begin{pmatrix}
        0 & 1 & 0 \\
        0 & 0 & 1 \\
        1 & 0 & 0 \\
        -1 & 0 & 0
    \end{pmatrix}.
\]
When $a_i$ and $b_i$ are both $0$-type, $\ell$ can be represented similarly.

Now assume that $a_1$ and $a_2$ are same type and so are $b_1$ and $b_2$. Replace $a_2-a_1$ by $\omega a_2'$ and $b_2-b_1$ by $\omega b_2'$. Then the representation of $\ell$ becomes
\[
    \ell =
    \begin{pmatrix}
        a_1 & a_2' & a_3 \\
        b_1 & b_2' & b_3
    \end{pmatrix}
    \terlattice
        2{\comega}0
        {\omega}30
        001
   \begin{pmatrix}
        \conj{a}_1 & \conj{b}_1 \\
        \conj{a}_2' & \conj{b}_2' \\
        \conj{a}_3 & \conj{b}_3
    \end{pmatrix}.
\]
This $3\times 3$ matrix can be represented as
\[
    \terlattice
        2{\comega}0
        {\omega}30
        001
    =
    \begin{pmatrix}
        0 & 0 & 0 & 1 \\
        0 & 1 & 1 & 0 \\
        1 & 0 & 0 & 0
    \end{pmatrix}
    \begin{pmatrix}
        1 & 0 & 0 & 0 \\
        0 & 1 & 0 & 0 \\
        0 & 0 & 2 & \omega \\
        0 & 0 & \comega & 2
    \end{pmatrix}
    \begin{pmatrix}
        0 & 0 & 1 \\
        0 & 1 & 0 \\
        0 & 1 & 0 \\
        1 & 0 & 0
    \end{pmatrix}.
\]
Thus if a pair of $a_i$ and $a_j$ is same type and a pair of $b_i$ and $b_j$ is same type, $\ell$ can be represented similarly. This argument can be applied to the case of that a pair of $a_i$ and $a_j$ have opposite types and a pair of $b_i$ and $b_j$ have opposite types.

Now let $F$ be the finite field $\frakO/\omega\frakO$. Consider $(a_i,b_i)$ as a vector in a vector space $F^2$ over $F$. We may assume that $(a_2,b_2)$ and $(a_3,b_3)$ are non-zero and linearly independent in $F^2$. Then we can find $\delta_1, \delta_2 \in \{1,-1\}$ such that $(a_3,b_3)-\delta_1(a_1,b_1)-\delta_2(a_2,b_2)$ is a pair of $0$-type numbers. Replace each component by $\omega a_3'$ and $\omega b_3'$. Then
\begin{align*}
    \ell &=
    \begin{pmatrix}
        a_1 & a_2 & \delta_1 a_1 + \delta_2 a_2 + \omega a_3' \\
        b_1 & b_2 & \delta_1 b_1 + \delta_2 b_2 + \omega b_3'
    \end{pmatrix}
    \terdiagonal111\\
    & \quad\quad\quad\quad\quad\quad\quad
    \begin{pmatrix}
        \conj{a}_1 & \conj{b}_1 \\
        \conj{a}_2 & \conj{b}_2 \\
        \delta_1 \conj{a}_1 + \delta_2 \conj{a}_2 + \comega
        \conj{a}_3'
        & \delta_1 \conj{b}_1 + \delta_2 \conj{b}_2 + \comega
        \conj{b}_3'
    \end{pmatrix}\\
    &=
    \begin{pmatrix}
        a_1 & a_2 & a_3' \\
        b_1 & b_2 & b_3'
    \end{pmatrix}
    \terlattice
        2{\delta_1\delta_2}{\delta_1\comega}
        {\delta_1\delta_2}2{\delta_2\comega}
        {\delta_1\omega}{\delta_2\omega}3
   \begin{pmatrix}
        \conj{a}_1 & \conj{b}_1 \\
        \conj{a}_2 & \conj{b}_2 \\
        \conj{a}_3' & \conj{b}_3'
    \end{pmatrix}.
\end{align*}

The $3 \times 3$ matrix can be represented as
\[
    \terlattice
        2{\delta_1\delta_2}{\delta_1\comega}
        {\delta_1\delta_2}2{\delta_2\comega}
        {\delta_1\omega}{\delta_2\omega}3
    =
    \begin{pmatrix}
        0 & 0 & 0 & \delta_1 \\
        0 & 0 & \delta_2 \comega & -\delta_2 \\
        0 & 1 & 1 & 0
    \end{pmatrix}
    \begin{pmatrix}
        1 & 0 & 0 & 0 \\
        0 & 1 & 0 & 0 \\
        0 & 0 & 2 & \omega \\
        0 & 0 & \comega & 2
    \end{pmatrix}
    \begin{pmatrix}
        0 & 0 & 0 \\
        0 & 0 & 1 \\
        0 & \delta_2 \omega & 1 \\
        \delta_1 & -\delta_2 & 0
    \end{pmatrix}.
\]

Hence the $2$-universality is proved.
\end{proof}


\section{Finiteness Theorem for 2-Universality}
\label{sec:FinitenessTheorem}

Recently, a beautiful criterion for universality of a given quadratic $\mathbb Z$-lattice was announced by Conway and Schneeberger \cite{Conway}. This criterion, known as the Fifteen Theorem, states\,: \emph{A positive definite quadratic $\mathbb Z$-lattice is universal if it represents every element in the set
\[
A:=\{\,1, 2, 3, 5, 6, 7, 10, 14, 15\,\}.%
\]}
Shortly after, an analogous criterion for $2$-universality was proved \cite{KKO1}, which states\,: \emph{A positive definite quadratic $\mathbb Z$-lattice is $2$-universal if it represents every element in the set
\[
B:=\left\{\,\qf{1,1}, \qf{2,3}, \qf{3,3}, \binlattice2112, \binlattice2113, \binlattice2114\,\right\}.%
\]}
We refer the readers to \cite{K} and \cite{KKO2} for recent developments in this direction.

The set of nine numbers is called \emph{minimal} in the sense that no proper subsets of those numbers ensure $2$-universality.

Now we investigate criteria as analogies of the 15-theorem. Criteria of this type are called \emph{finiteness theorems}.

\begin{Thm}
    If a Hermitian lattice over $\Q{-1}$ represents $\left\{ \qf{1,1},     \binlattice{2}{1}{1}{2} \right\}$, it is $2$-universal.
\end{Thm}
\begin{proof}
Let $L$ be a $2$-universal lattice over $\Q{-1}$. Let $\{\mathbf{v}_1,\mathbf{v}_2 \}$ and $\{ \mathbf{v}_3, \mathbf{v}_4 \}$ be the bases of $\qf{1,1}$ and $\binlattice2112$, respectively. Then $L$ contains a lattice generated by all $\mathbf{v}_i$'s. This lattice can be obtained by using the following positive semi-definite $4\times4$-matrix
\[
    \begin{pmatrix}
        1 & 0 & \ast & \ast \\
        0 & 1 & \ast & \ast \\
        \ast & \ast & 2 & 1 \\
        \ast & \ast & 1 & 2
    \end{pmatrix}.
\]
It is isometric to $\qf{1,1,1,0}$, $\qf{1,1,1,1}$, or $\qf{1,1}\perp\binlattice{2}{1}{1}{2}$. Thus $L$ contains $\qf{1,1,1}$ or $\qf{1,1}\perp\binlattice{2}{1}{1}{2}$. Both lattices are $2$-universal.
\end{proof}

\begin{Thm}
    If a Hermitian lattice over $\Q{-2}$ represents $\left\{     \qf{1,1}, \binlattice{2}{-1+\omega}{-1+\comega}2 \right\}$, it is     $2$-universal.
\end{Thm}
\begin{proof}
The positive semi-definite matrix
\[
    \begin{pmatrix}
        1 & 0 & \ast & \ast \\
        0 & 1 & \ast & \ast \\
        \ast & \ast & 2 & -1+\omega \\
        \ast & \ast & -1+\comega & 2
    \end{pmatrix}
\]
gives only one lattice $\qf{1,1}\perp\binlattice{2}{-1+\omega}{-1+\comega}2$. This lattice was proved to be $2$-universal.
\end{proof}

\begin{Thm}\label{Thm:finiteness}
    Let $E=\Q{-3}$ and $\frakO$ be its ring of integers. A Hermitian     $\frakO$-lattice is $2$-universal if it represents $\qf{1,1}$ and     $\qf{1,2}$.
\end{Thm}

\begin{proof}
Let $L$ be a Hermitian $\frakO$-lattice and assume that $\qf{1,1}\rightarrow L$ and $\qf{1,2}\rightarrow L$. Since $\qf{1,1}$ is unimodular, it splits $L$, that is, $L=\qf{1,1}\perp L_0$ for some sublattice $L_0$ of $L$. In order for $L$ to represent $\qf{1,2}$, $L_0$ should represent $1$, $2$, or $\qf{1,2}$. In any case, $L$ contains either $\qf{1,1,1}$ or $\qf{1,1,2}$. Since both are $2$-universal, so is $L$.
\end{proof}

\begin{Rmk}
The sets $A$ and $B$ above are unique minimal sets in the respective criteria. In Theorem \ref{Thm:finiteness}, however, $\qf{1,2}$ can be replaced by
\[
    \binlattice2{\omega}{\comega}3.
\]
So, the set $\{\,\qf{1,1},\qf{1,2}\,\}$ is a minimal but not a unique set ensuring the $2$-universality of $L$.
\end{Rmk}

\begin{Thm}
    If a Hermitian lattice over $\Q{-11}$ represents $\left\{     \qf{1,1}, \binlattice{2}{\omega}{\comega}2 \right\}$, it is     $2$-universal.
\end{Thm}
\begin{proof}
The positive semi-definite matrix
\[
    \begin{pmatrix}
        1 & 0 & \ast & \ast \\
        0 & 1 & \ast & \ast \\
        \ast & \ast & 2 & \omega \\
        \ast & \ast & \comega & 2
    \end{pmatrix}
\]
gives only one lattice
$\qf{1,1}\perp\binlattice{2}{\omega}{\comega}2$. This lattice was proved to be $2$-universal.
\end{proof}

In summary we conclude that the following are all \emph{new} $2$-universal Hermitian lattices over $\Q{-m}$ when $m=1,2,3,11$. 

\begin{center}
\begin{tabular}{ll}
$\Q{-1}$\,: & $\qf{1,1,1}$, $\qf{1,1}\perp\binlattice{2}{1}{1}{2}$ \\
$\Q{-2}$\,: & $\qf{1,1}\perp\binlattice{2}{-1+\omega_2}{-1+\comega_2}{2}$ \\
$\Q{-3}$\,: & $\qf{1,1,1}$, $\qf{1,1,2}$ \\
$\Q{-11}$\,: & $\qf{1,1} \perp
\binlattice{2}{\omega_{11}}{\comega_{11}}{2}$
\end{tabular}
\end{center}

\section{2-Universal Hermitian Lattices of Higher Rank}
We denote the minimal rank of $2$-universal Hermitian lattices over $\Q{-m}$ by $u_2(-m)$. We know that $u_2(-1) = 3$, $u_2(-2) = 4$, $u_2(-3) = 3$, $u_2(-7) = 3$, and $u_2(-11) = 4$.

Now assume $m=5,6,10$ or $m\ge13$ in this section. It is clear that there exist $2$-universal Hermitian lattices for all $m$, because the lattice $I_3=\qf{1,1,1}$ is locally $2$-universal and thus we can make a $2$-universal Hermitian lattice by summing up orthogonally all classes in the genus of $I_3$. We obtain an upper bound for minimal rank as follows:
\[
    u_2(-m) \le 3 h(I_3).
\]

For example, from \cite{Schiemann} we have that over $\Q{-19}$
\[
    \gen I_3 = \left\{ \terdiagonal111, \terlattice10002\omega0\comega3,
    \terlattice21113{1+\omega}1{1+\comega}3 \right\}.
\]
Thus the following lattice of rank 8 is $2$-universal over
$\Q{-19}$.
\[
    \qf{1,1,1} \perp \binlattice{2}{\omega}{\comega}{3} \perp
    \terlattice{2}{1}{1} {1}{3}{1+\omega} {1}{1+\comega}{3}
\]
%

Let us consider a diagonal Hermitian form $L = \qf{a_1,\dotsc,a_n}$. Choose an integer $c > 2$ satisfying $\frac{\omega\comega}2 < c < \omega\comega$. This is possible since $\omega\comega \ge 4$. We construct a target lattice $\ell = \binlattice{2}{\omega}{\comega}{c}$, which is positive definite. If $L$ represents $\ell$, then the representation has a term $a_i(\alpha x+\beta y)\conj{(\alpha x+\beta y)}$ with $\im\alpha\ne0$ or $\im\beta\ne0$ for some $i$ because of the cross term $\omega x\conj{y}$. But then
\[
    (\alpha x+\beta y)\conj{(\alpha x+\beta y)}
     = \alpha\conj{\alpha} x\conj{x} +
    \alpha\conj{\beta} x\conj{y} + \conj{\alpha}\beta \conj{x}y +
    \beta\conj{\beta}y\conj{y}
\]
and $\alpha\conj{\alpha} > c$ or $\beta\conj{\beta} > c$. Thus $L$ cannot represents $\ell$.

\begin{Rmk}
Unlike quadratic forms over $\Z$, there are no \emph{diagonal} $2$-universal Hermitian forms when $m=5,6,10$ or $m\ge13$.
\end{Rmk}


Now consider binary sublattices $\ell(k) := \binlattice{k}\omega\comega{c_k}$ with $2 \le k \le c_k$ and $\frac{\omega\comega}k < c_k \le \frac{\omega\comega}{k-1}$. Then the largest index of $k$'s is given as $n := \left\lfloor \frac{\sqrt{4\omega\comega+1}+1}2 \right\rfloor$ because
\[
    n \le \frac{\omega\comega}{n-1} \text{ and }
    \frac{\omega\comega}n+1 \le \frac{\omega\comega}{n-1}.
\]
These become
\[
    2 \le n \le \frac{1+\sqrt{1+4\omega\comega}}2.
\]
 
Let $\{\mathbf{v}_k, \mathbf{w}_k\}$ be the basis of $\ell(k)$ with $\mathbf{v}_k \cdot \mathbf{v}_k = k$, $\mathbf{w}_k \cdot \mathbf{w}_k = c_k$, and $\mathbf{v}_k \cdot \mathbf{w}_k = \omega$. Also let $\mathbf{v}_1 \cdot \mathbf{v}_1 = 1$. We show that $\mathbf{v}_1, \cdots, \mathbf{v}_n$ are linearly independent.

Suppose that those vectors are not linearly independent. Then
\[
a_1 \mathbf{v}_1 + \cdots + a_n \mathbf{v}_n + \omega(b_1 \mathbf{v}_1 + \cdots + b_n \mathbf{v}_n) = 0
\]
with $a_i, b_i \in \Z$. Note that $|\mathbf{v}_i \cdot \mathbf{v}_j|^2 \leq ij < \omega\comega$ and thus $\mathbf{v}_i \cdot \mathbf{v}_j \in \Z$. Multiplying by $\mathbf{v}_j$ and comparing both sides, we conclude that $(b_1 \mathbf{v}_1 + \cdots + b_n \mathbf{v}_n)\cdot\mathbf{v}_j = 0$. Since the norm of $b_1 \mathbf{v}_1 + \cdots + b_n \mathbf{v}_n$ vanishes and the concerned $2$-universal Hermitian lattice is positive definite, $b_1
\mathbf{v}_1 + \cdots + b_n \mathbf{v}_n = 0$. Thus we can write
\[
    a_{k_1} \mathbf{v}_{k_1} + \cdots + a_{k_N} \mathbf{v}_{k_N} = 0
\]
with nonzero $a_{k_i} \in \Z$ and $k_1 < k_2 < \cdots < k_N$. But, we obtain a contradiction by multiplying both sides by $\mathbf{w}_{k_N}$, since $\mathbf{v}_{k_i} \cdot \mathbf{w}_{k_N} \in \Z$ for $i < N$ and $\mathbf{v}_{k_N} \cdot \mathbf{w}_{k_N} = \omega$. So $\mathbf{v}_1, \cdots, \mathbf{v}_n$ are linearly independent and ranks of $2$-universal lattices are bigger than $n$. That is, $u_2(-m) > \left\lfloor \frac{\sqrt{4\omega\comega+1}+1}2 \right\rfloor$.

If we denote the discriminant of $E=\Q{-m}$ by $d_E$, $d_E = 4m$ if $m \nequiv 3\pmod{4}$ and $d_E = m$ if $m \equiv 3 \pmod{4}$. Then the above inequality implies
\[
    u_2(-m) = u_2(d_E) = \Omega(\sqrt{d_E}),
\]
where $f(n) = \Omega(g(n))$ means $\lim\inf_{n\to\infty} \mid{f(n)}/{g(n)}\mid > 0$. Roughly speaking, $u_2(d_E)$ increases as $d_E$ increases. In addition, so does $h(I_3)$.

Let us think about finiteness theorems for $2$-universal lattices of higher ranks. We have constructed binary lattices $\binlattice{k}{\omega}{\comega}{c_k}$ which $L$ must represent. Thus if there exists a set $S_{-m}$ ensuring $2$-universality over $E=\Q{-m}$, the cardinality of $S_{-m}$ has a lower bound:
\begin{align*}
    \#S_{-m} & \ge \text{the number of $c_n$'s}+2 \\
    &=
    \begin{cases}\displaystyle
        \left\lfloor \frac{\sqrt{4m+1}-1}2 \right\rfloor + 1 & \text{if
        } m \nequiv 3 \pmod{4}, \\ \displaystyle
        \left\lfloor \frac{\sqrt{m+2}-1}2 \right\rfloor + 1 & \text{if } m
        \equiv 3 \pmod{4}.
    \end{cases}
\end{align*}
In other words, $\#S_{d_E} = \Omega(\sqrt{d_E})$. The cardinality of $S_{d_E}$ increases as $d_E$ increases.


%

\end{document}